\documentclass[a4paper,10pt]{smfart}
\usepackage[french]{babel}
\usepackage{t1enc}
\usepackage{a4}
\usepackage{amsfonts,amssymb}
\usepackage{amsmath}
\usepackage{amsopn}
\usepackage{amsthm}

\usepackage{latexsym}
\usepackage{bull}
\newcommand {\be }{ \begin{equation} }
\newcommand {\ee }{ \end{equation} }
\newcommand {\ben }{ \begin{equation*} }

\newcommand {\een }{ \end{equation*} }
\newcommand {\bea }{ \begin{eqnarray} }
\newcommand {\eea }{ \end{eqnarray} }
\newcommand {\ba }{\begin{array}}
\newcommand {\ea} {\end{array}}

\newcommand {\lbl }{\label}
\newcommand {\ul }{\underline}

\newcommand {\C }{\mathbb{C}}
\newcommand {\N }{\mathbb{N}}
\newcommand {\Z }{\mathbb{Z}}

\newcommand {\ml } {\left( \begin{array}}
\newcommand {\mr} {\end{array} \right) }
\newcommand {\ts} {\widetilde{\sigma _{\lambda}}}
\newcommand {\Us} {_{\sigma _{\rho}} \! U_q(\mathfrak g)  }
\newcommand {\sla} {\sigma _{\lambda}}
\newcommand {\sr} {\sigma _{\rho}}\newcommand {\ds} {\cdot _{\sigma _{\rho}}}
\newcommand {\Ug} {U_q(\mathfrak g )}
\newcommand {\grU} {\mathrm{Gr}\; U_q( \mathfrak g )} \newcommand {\Gal} {\mathrm{Gal}}
\newcommand {\cH} {\mathcal{H}}

\newtheorem{theorem}{Th\'eor\`eme}
\newtheorem*{theo}{Th\'eor\`eme}
\newtheorem{lemma}[theorem]{Lemme}

\title{Classification des objets galoisiens de $U_q(\mathfrak g)$
à homotopie près}
\alttitle{Classification of Galois objects of~$U_q(\mathfrak g)$ up to homotopy equivalence}
\date{10.10.2006}

\author{Thomas AUBRIOT}
\address{Institut de Recherche Mathématique Avancée \\C.N.R.S. - Université Louis Pasteur
\\7 rue René Descartes \\67084 Strasbourg Cedex, France \\Fax : +33 (0)3 90 24 03 28}
\email{aubriot@math.u-strasbg.fr}

\begin{document}
\setlength{\textwidth}{13cm}
\setlength{\textheight}{21.5cm}

\begin{abstract}
Pour toute algèbre enveloppante quantique $\Ug$ de Drinfeld-Jimbo
et toute famille~$\lambda =(\lambda _{ij})_{1\leq i<j\leq t} \in
k^{\star}$ d'éléments inversibles du corps de base, nous construisons
explicitement par générateurs et relations un objet
galoisien~$A_{\lambda}$ de~$\Ug$ et nous montrons que tout objet galoisien
de~$\Ug$ est homotope à un unique objet de la forme~$A_{\lambda}$.
\end{abstract}

\begin{altabstract} For any Drinfeld-Jimbo quantum enveloping algebra~$\Ug$ and
for any family~$\lambda =(\lambda _{ij})_{1\leq i<j\leq t} \in
k^{\star}$ of invertible elements of the base field, we explicitly
construct a Galois object~$A_{\lambda}$ of~$\Ug$ by generators and relations and we prove that any Galois object of~$\Ug$
is homotopic to a unique object of type~$A_{\lambda}$.
\end{altabstract}

\altkeywords{Galois extension, Hopf algebra, Drinfeld-Jimbo quantum group, Homotopy, Noncommutative geometry, Principal fibre bundle}
\maketitle \mainmatter

\section*{Introduction}

Le concept d'extension Hopf-galoisienne qui a \'et\'e beaucoup
\'etudi\'e ces derni\`eres ann\'ees
est une g\'en\'eralisation naturelle du concept classique d'extension
galoisienne de corps commutatifs.
C'est aussi l'analogue alg\'ebrique de la notion de fibr\'e principal
dans le cadre de la g\'eom\'etrie
non commutative.

Bien qu'une litt\'erature abondante ait \'et\'e consacr\'ee aux
extensions Hopf-galoi\-siennes
(voir par exemple \cite{M}, \cite{S} et les r\'ef\'erences donn\'ees dans ces
deux articles),
on a peu de r\'esultats sur leur classification \`a isomorphisme pr\`es.
Pour contourner la difficult\'e --- qui semble grande ---
de classer les extensions Hopf-galoisiennes \`a isomorphisme pr\`es,
Kassel \cite{K} a introduit sur
les extensions Hopf-galoisiennes une relation d'équivalence moins fine que
l'isomorphie, relation
qu'il a appel\'ee homotopie.
Dans \cite{KS} Kassel et Schneider en ont fait une \'etude syst\'ematique.
Ils donnent notamment
une application \`a la classification des extensions Hopf-galoisiennes
lorsque l'alg\`ebre de Hopf
est l'alg\`ebre enveloppante quantique $\Ug$ associ\'ee par Drinfeld
et Jimbo
\`a une alg\`ebre de Lie semi-simple complexe $\mathfrak g$.
L'une des cons\'equences des r\'esultats de \cite{KS}
porte sur l'ensemble~$\cH_k(\Ug)$ des classes d'homotopie
des extensions~$\Ug$-galoisiennes du corps de base~$k$
(ces objets sont \'egalement appel\'es objets galoisiens de~$\Ug$) :
Kassel et Schneider d\'emontrent que $\cH_k(\Ug)$ est en bijection
avec
le groupe de cohomologie $H^2(G,k^*)$,
o\`u $G$ est le groupe des \'el\'ements ``group-like'' de~$\Ug$ opérant
trivialement sur le groupe~$k^{\star}$ des éléments inversibles de~$k$.
Le groupe~$G$ est un groupe ab\'elien libre dont le rang $t$ est \'egal
\`a celui
de l'alg\`ebre de Lie~$\mathfrak g$.
Il est bien connu que tout \'el\'ement de~$H^2(G,k^*)$ peut \^etre
repr\'esenté
par une famille~$\lambda$ de~$t(t-1)/2$ \'el\'ements non nuls du
corps de base.

Le but de cet article est de construire explicitement par
g\'en\'erateurs et relations
un objet galoisien~$A_{\lambda}$ de~$\Ug$ pour toute famille~$\lambda$ de ce type
et de montrer que tout objet galoisien de~$\Ug$ est homotope \`a un
unique
objet galoisien de la forme~$A_{\lambda}$.

Au paragraphe $1$, nous rappelons la d\'efinition des concepts
d'extension Hopf-galoisienne et d'objet galoisien. Nous redonnons
\'egalement la pr\'esentation standard de l'alg\`ebre enveloppante
quantique~$\Ug$.

Les objets galoisiens~$A_{\lambda}$ sont construits au paragraphe $2$.
Nous y \'enon\c cons aussi le th\'eor\`eme principal de l'article.

Le paragraphe $3$ est enti\`erement consacr\'e \`a la d\'emonstration du
th\'eor\`eme.

\section{Rappels}

\subsection{Extensions galoisiennes et objets galoisiens}

Soit $k$ un corps commutatif. Tous les objets de cet article
appartiennent à la cat\'egorie tensorielle des~$k$-espaces
vectoriels et nous ne considérons que des algèbres de Hopf
admettant une antipode bijective.
Si $H$ est une alg\`ebre de Hopf et $A$ est une
alg\`ebre $H$-comodule \`a droite dont la coaction est le
morphisme d'alg\`ebres~$\delta : A \to A\otimes H$, nous
d\'efinissons la sous-alg\`ebre $B$ des \'el\'ements
$H$-covariants de~$A$ par
\be B = \{ a\in A \mid \delta(a) =a\otimes 1\} .\end{equation}
L'application lin\'eaire $\beta : A \otimes_B A
\to A \otimes H$ d\'efinie par
\be \beta(a\otimes a') = (a\otimes 1)\delta(a'), \end{equation}
 pour $a$, $a'\in A$, est appel\'ee {\it l'application canonique} associ\'ee \`a~$A$.
Une alg\`ebre~$H$-comodule \`a droite~$A$ est une {\it extension
$H$-galoisienne}
de~$B$ si~$B$ est la sous-alg\`ebre des \'el\'ements $H$-covariants~de
$A$, si l'application canonique~$\beta : A \otimes_B A \to H \otimes A$
associ\'ee \`a~$A$
est un isomorphisme et si~$A$ est fid\`element plat en tant que~$B$-module \`a droite ou \`a gauche.

Un {\it objet galoisien} d'une alg\`ebre de Hopf $H$ est une
extension $H$-galoisienne
du corps de base~$k$.

Deux extensions $H$-galoisiennes $A$ et $A'$ de $B$ sont dites {\it isomorphes}
s'il existe
un morphisme $f: A \to A'$ d'alg\`ebres $H$-comodules qui soit un
isomorphisme
et qui soit l'identit\'e sur~$B$.

Nous notons $\Gal _B(H)$ l'ensemble des classes d'isomorphisme d'extensions $H$-galoisiennes de $B$. L'ensemble $\Gal _B(H)$ peut
\^etre consid\'er\'e comme un foncteur contravariant en $H$. En
effet, soit $i : K \to H$ un morphisme d'alg\`ebres de Hopf.
Rappelons \cite{EM} que, étant donné une algèbre de Hopf $H$, un
comodule $A$ à droite de coaction $\delta _A$ et un comodule $K$ à
gauche de coaction $\delta _K$, le produit cotensoriel~$A\Box _{H}
K$ est défini comme le noyau de l'application \be \mathrm{Id}_A \otimes
 \delta_K - \delta _A \otimes \mathrm{Id}_K: A\otimes K
\rightarrow A\otimes H \otimes K ,\end{equation}
 (ou encore l'égalisateur des coactions de $A$ et
$K$). Si~$A$ est une
extension $H$-galoisienne de~$B$ à droite, alors \be i^{\star}
(A) = A\Box _H K \end{equation} est une extension $K$-galoisienne de $B$ à droite d'après \cite[Prop 3.11
(3)]{S}.

Kassel et Schneider \cite{KS} (voir aussi \cite{K}) ont d\'efini
une relation d'\'equivalence, not\'ee $\sim$ et appel\'ee {\it
homotopie}, sur l'ensemble $\Gal_B(H)$ (nous renvoyons \`a
\cite{KS} pour la d\'efinition). Nous notons~$\cH_B(H)$ l'ensemble
des classes d'homotopie d'extensions~$H$-galoisiennes de~$B$.
L'application \be i^* : \Gal_B(H) \to \Gal_B(K)\end{equation} induite par un
morphisme d'alg\`ebres de Hopf $i : K \to H$ et d\'efinie plus
haut, passe aux classes d'homotopie et d\'efinit une application
\be i^* : \cH_B(H) \to \cH_B(K). \end{equation}

\subsection{Cocycles et extensions clivées}
Suivant (\cite[Chapitre 7]{M}), nous dirons qu'une application
bilinéaire $\sigma : H\times H \rightarrow k$ est un {\it cocycle
normalisé inversible} pour l'algèbre de Hopf $H$ si $\sigma$ est inversible pour la convolution et vérifie les
relations \be \lbl{cocycle} \sigma  (x_{(1)},y_{(1)}) \sigma
(x_{(2)}y_{(2)},z)= \sigma  (y_{(1)},z_{(1)}) \sigma
(x,y_{(2)}z_{(2)}) \end{equation} et \be \lbl{norm} \sigma  (1,x)= \sigma
(x,1)=\varepsilon (x),\end{equation} pour $x,y,z \in H$ ($\varepsilon $ est la coünité de $H$).

Nous avons utilisé ici la notation de Sweedler $ \Delta
(x)=x_{(1)} \otimes x_{(2)} $ pour la comultiplication $\Delta$ de
$H$, notation que nous utiliserons dans la suite de l'article pour les
différentes comultiplications et coactions.

Rappelons (\cite[Chapitre 7]{M}) que si $H$ est une algèbre de Hopf,
$\sigma :H\times H \rightarrow k $ un cocycle inversible normalisé et
$B$ une algèbre, alors le {\it produit croisé}~$ B \sharp _{\sigma} H$ est
l'espace vectoriel $B\otimes H$ muni du produit associatif et unifère
\be
\lbl{produit_croise} (a\sharp h)(b\sharp k) = \sigma (h_{(1)}
,k_{(1)}) ab \sharp h_{(2)} k_{(2)} ,\end{equation}
pour tout $a,b \in B$ et $h,k \in H$.
De plus, cette algèbre peut \^etre munie d'une structure d'algèbre $H$-comodule à droite induite
par la comultiplication de $H$, ce qui fait de $B\sharp _{\sigma } H$ une extension~$H$-galoisienne de
$B$. Les extensions $H$-galoisiennes de cette forme sont appelées
{\it extensions clivées} ({\it cleft} en anglais). Lorsque $k=B$, le produit croisé~$k\sharp _{\sigma} H$
s'identifie à $H$ muni du produit \be
\lbl{produit_tordu_clive_gauche} x \cdot _{\sigma} y = \sigma (x_{(1)} ,y_{(1)} ) x_{(2)}
y_{(2)}, \end{equation} pour $x,y \in H$. Nous notons
$_{\sigma}H$ l'espace vectoriel~$H$ muni de ce produit associatif dont l'unité est celle de $H$.

Si $H$ est une algèbre de Hopf et $\rho$ est un cocycle normalisé inversible, nous pouvons aussi définir l'algèbre de Hopf~$H^{\rho}$ comme la
cogèbre $H$ munie du produit associatif modifié
\be
\lbl{produit_tordu_hopf} x \cdot _{\rho} y = \rho (x_{(1)} ,y_{(1)})  x_{(2)}
y_{(2)} \rho ^{-1} (x_{(3)} ,y_{(3)}), \end{equation} pour $x,y \in H$. Nous définissons aussi pour toute
$H$-comodule algèbre à droite $A$ la $H^{\rho}$-comodule algèbre tordue à droite $A_{\rho}$ qui est le $H^{\rho}$-comodule (ou $H$-comodule)
$A$ muni du produit
\be \lbl{produit_tordu_comod_droite} a \cdot _{\rho} b = a_{(1)} b_{(1)}  \rho ^{-1}(a_{(2)}, b_{(2)}), \end{equation}
pour $a,b \in A$. Alors il existe une bijection entre les extensions~$H$-galoisiennes de $B$ à droite et les extensions $H^{\rho}$-galoisiennes de $B$  à droite~$j: \Gal _B (H) \to \Gal _B (H^{\rho})$ donnée par \be j(A) =  A_{\rho}. \end{equation}

Montgomery et Schneider \cite[Théorème 5.3]{MS} ont montré que, si $_{\sigma}H$ est une extension clivée de $H$, alors son image
$j(_{\sigma}H)=(_{\sigma}H)_{\rho}$ dans $\Gal _k (H^{\rho})$ est une extension clivée  de~$H^{\rho}$ qui s'identifie à~$_{\sigma *\rho ^{-1}}(H^{\rho})$.

Au paragraphe 3, nous aurons besoin du lemme suivant. 
\begin{lemma}
\lbl{lemme_iUs}
Soit $H'$ une sous-algèbre de Hopf d'une algèbre de Hopf $H$ et $\sigma ' : H'\otimes H'\to k$ un cocycle de $H'$, s'étendant en un cocycle
$\sigma  : H\otimes H\to k$ de $H$. Nous avons alors un isomorphisme d'algèbres
\begin{equation}
  _{\sigma } H \Box _H H' \cong \; _{\sigma '} H' .
\end{equation}
\end{lemma}

\begin{proof}
Considérons les bijections linéaires (voir \cite[Prop 2.1]{EM})
\be H' \xrightarrow{\mu } H \Box _H H' \xrightarrow{\nu} H'  \end{equation} données par
$\mu = (i\otimes \mathrm{Id}) \circ \Delta _{H'} $ et $\nu =\varepsilon
\otimes \mathrm{Id}$. Les structures d'espaces vectoriels et de
comodules de~$_{\sigma }H $ et $_{\sigma '} H'$ étant par définition les
m\^emes que celles de~$H$ et $H'$ respectivement, ces
bijections valent aussi pour~$_{\sigma }H$ et~$_{\sigma '} H'$ : \be
_{\sigma '} H' \xrightarrow{\mu} \; _{\sigma }H \Box _{H}H'
\xrightarrow{\nu} \: _{\sigma '} H' .  \end{equation}

Montrons maintenant que l'isomorphisme $\mu$ est aussi un
morphisme d'algèbres de $_{\sigma '} H'$ sur $_{\sigma} H \Box _{H} H$.
Notons $g \cdot _{\sigma ' }h$ le produit dans $_{\sigma '}\! H'$ de
deux éléments $g$ et $h$ de $H'$ ; nous gardons la notation
$gh$ pour leur produit dans $H'$ et les notations $g,h$ pour les éléments $g,h\in H'$ vu dans $H$. Nous notons aussi
$g \cdot _{\sigma  }h$ le produit dans $_{\sigma }\! H$ de deux éléments $g,h$.
Nous obtenons
\be \ba{rcl} \mu (g \cdot _{\sigma ' } h )&=& \sigma ' (h_{(1)},g_{(1)}) \mu (g_{(2)}h_{(2)}) \\
&=& \sigma '(h_{(1)},g_{(1)}) ( g_{(2)} h_{(2)} \otimes g_{(3)}h_{(3)}) \\
&=& (\sigma (h_{(1)},g_{(1)})  g_{(2)} h_{(2)}) \otimes g_{(3)}h_{(3)} \\
&=& g_{(1)} \cdot _{\sigma  }  h_{(1)} \otimes g_{(2)} h_{(2)} \\
\ea \end{equation} pour tout $ g,h \in G$, ce qui assure que
\be \mu : \: _{\sigma '}\! H' \to  \; _{\sigma }H \Box _H H' \end{equation} est un
isomorphisme d'algèbres.

\end{proof}

\subsection{Les alg\`ebres enveloppantes quantiques de Drinfeld-Jimbo}

Nous supposons désormais que $k$ est un corps de caractéristique
différente de~$2$ ou~$3$. Fixons la matrice de Cartan~$(a_{ij})_{1\leq i,j \leq t}$ d'une algèbre de Lie
semi-simple complexe $\mathfrak g$, des entiers~$(d_i)_{1\leq
i\leq t} \in \{1,2,3\}$ tels que~$d_i a_{ij}= d_j a_{ji}$ pour
tout~$1\leq i,j\leq t$ ainsi qu'un élément inversible~$q \in
k $ tel que $q^{2d_i}\neq 1$ pour tout~$i=1,\ldots , t$.

L'algèbre de Drinfeld-Jimbo $\Ug$ (voir \cite[Chapitre 4]{J}) est l'algèbre associative unitaire
engendrée par les générateurs $E_i,F_i,K_i$ et $K_i^{-1}$ pour $1\leq i \leq t$ et les
relations

\be K_i  K_j = K_j  K_i ,
\quad
K_i  K_i^{-1} = K_i ^{-1}  K_i =1  ,
\lbl{comK_K}
\end{equation}

\be K_i  E_j = q^{d_ia_{ij}} E_j K_i , \lbl{comK_E}\end{equation}

\be K_i  F_j = q^{-d_i a_{ij} } F_j K_i,
 \lbl{comK_F}\end{equation}

\be E_i  F_j -F_j  E_i = \delta _{ij}
\frac{K_i-K_i^{-1}}{q^{d_i}-q^{-d_i}}
\lbl{comE_F},\end{equation}

\be \sum _{r=0}^{1-a_{ij}} (-1)^r \left[ \ba{c} 1-a_{ij} \\ r \ea
\right]_{q^{d_i}} E_i^{1-a_{ij}-r}
 E_j E_i ^r =0 ,
\lbl{SerreE}
\end{equation}

\be \sum _{r=0}^{1-a_{ij}} (-1)^r \left[ \ba{c} 1-a_{ij} \\ r \ea
\right]_{q^{d_i}} F_i^{1-a_{ij}-r} F_j F_i ^r =0,
 \lbl{SerreF} \end{equation}
pour $1\leq i,j\leq t $.

Il est bien connu (\cite[Chapitre 4]{J}) que $U_q(\mathfrak g)$
peut \^etre munie  d'une
structure d'algèbre de Hopf
avec la comultiplication $\Delta $ définie sur les
générateurs par
\be \lbl{comul_U_qE}
\Delta (E_i)=E_i \otimes 1 + K_i \otimes E_i ,
 \end{equation}
\be \lbl{comul_U_qF}
\Delta  (F_i)= F_i \otimes K_i^{-1} + 1 \otimes F_i,
\end{equation}
\be \lbl{comul_U_qK}
\Delta  (K_i ^{\pm 1})=  K_i ^{\pm 1} \otimes K_i ^{\pm 1} , \end{equation}
la co\"unité $\varepsilon $ définie par
\be \lbl{counit_U_q}
\varepsilon (E_i)=0,\quad
\varepsilon  (F_i)=0 ,\quad
\varepsilon (K_i ^{\pm 1})= 1 ,
\end{equation}
et l'antipode $S$ définie par
\be \lbl{anti_U_q}
S (E_i)=-K_i^{-1} E_i,\quad
S  (F_i)= -F_iK_i ,\quad
S (K_i ^{\pm 1})= K_i^{\mp 1}
,\end{equation}
pour tout $1\leq i \leq t$.

Notons $G$ le sous-groupe multiplicatif de~$U_q(\mathfrak g)$
engendré par~$K_1,K_2,\ldots ,K_t$. C'est un sous-groupe abélien libre de rang $t$. Notons encore $\Ug ^+$
la sous-algèbre de~$\Ug$ engendrée par les générateurs $K_i, K_i^{-1}$ et $E_i$, pour $1 \leq i \leq t$, et les relations
(\ref{comK_K}), (\ref{comK_E}) et (\ref{SerreE}) et~$\Ug ^-$ celle engendrée par les générateurs~$K_i,K_i^{-1}$ et $F_i$, pour $1 \leq i \leq t$, et les relations (\ref{comK_K}), (\ref{comK_F}) et (\ref{SerreF}).
L'algèbre de Hopf $\Ug$ est filtrée avec les générateurs $K_i^{\pm 1}$ en degré~$0$ et les générateurs~$E_i, F_i$ ~en degré~$1$. L'algèbre de Hopf graduée $\grU$ associée à cette filtration est engendrée, comme algèbre, par $E_i,F_i,K_i$
et~$K_i^{-1}$, pour $1\leq i \leq t$, soumis aux relations (\ref{comK_K}) - (\ref{comK_F}), (\ref{SerreE}), (\ref{SerreF})
ainsi qu'à la relation de
commutation \be \lbl{comE_F_grU} E_i F_j -F_j E_i =0,\end{equation} pour $1\leq i,j \leq t$.
Notons aussi $\grU ^+$
la sous-algèbre de $\grU$ engendrée par les générateurs $K_i, K_i^{-1}$ et $E_i$, pour $1 \leq i \leq t$, et les relations
(\ref{comK_K}), (\ref{comK_E}) et (\ref{SerreE}) et~$\grU ^-$ celle engendrée par les générateurs
$K_i,K_i^{-1}$ et $F_i$, pour~$1 \leq i \leq t$, et les relations~(\ref{comK_K}), (\ref{comK_F}) et~(\ref{SerreF}) et
remarquons que nous pouvons identifier $\grU ^+ \cong \Ug ^+$ ainsi que $\grU ^- \cong \Ug ^-$.

Kassel et Schneider \cite[Section 4]{KS} ont montré qu'il existe un unique cocycle normalisé $\rho : \grU \times \grU \to k$ qui vérifie
\be \lbl{rhoK_K} \rho (K_i,K_j)= 1, \end{equation} \be \lbl{rhoE_F} \rho (E_i ,F_j)= -\frac{\delta _{ij}}{q^{d_i}-q^{-d_i}},
\end{equation} pour $1\leq i,j \leq t$ et 
\be \lbl{rho_gene} \rho (x,y) =0, \end{equation} pour tous les autres couples de générateurs $x,y$ (en particulier ce cocycle n'est pas symétrique et $\rho (F_j, E_i)=0$). De plus si $x,y,z$ appartiennent à la même sous-algèbre $\grU ^+$ ou $\grU ^-$, on a
\be \lbl{multiplication1_rho}\rho (x,yz) = \rho (x_{(1)},y) \rho (x_{(2)},z), \end{equation}
\be \lbl{multiplication2_rho}\rho (xy,z) = \rho (x,z_{(2)}) \rho (y,z_{(1)}).\end{equation}
La relation (\ref{multiplication1_rho}) est encore vraie si $x\in \grU ^-$ et $y,z \in \grU ^+$ et la relation (\ref{multiplication2_rho}) est vraie si $x,y \in \grU ^-$ et $z\in \grU ^+$.
En particulier, ces relations impliquent
\be \rho (x,1)=\rho (1,x) =\varepsilon (x), \end{equation}
pour tout $x\in \grU$.
Kassel et Schneider \cite[Section 4]{KS} ont établi qu'il existe un isomorphisme d'algèbres de Hopf $\Ug \cong \grU ^{\rho}$ qui est l'identité sur les générateurs.

\section{Le résultat}

Nous considérons une famille~$(\lambda _{ij})_{1\leq i <j \leq t}$
d'éléments inversibles de~$k$. Par commodité, nous posons~$\lambda
_{ij}=\lambda _{ji}^{-1}$ et~$\lambda _{ii} =1 $ pour tout~$1 \leq
j\leq i \leq t$.

Nous définissons l'algèbre $A_{\lambda}$ comme l'algèbre
associative unitaire
engendrée par des générateurs $X_i ,Y_i,Z_i,Z_i^{-1}$ pour $1\leq i\leq t$
et les relations

\be Z_i  Z_j = \lambda _{ij}^2 Z_j  Z_i , \quad Z_i  Z_i^{-1} =
Z_i ^{-1}  Z_i =1  , \lbl{comZ_Z} \end{equation}

\be Z_i  X_j =\lambda _{ij}^{2} q^{d_ia_{ij}} X_j Z_i ,
\lbl{comZ_X}\end{equation}

\be Z_i  Y_j = q^{-d_i a_{ij} } Y_j  Z_i
 \lbl{comZ_Y},\end{equation}

\be X_i  Y_j -Y_j  X_i = \delta _{ij}
\frac{Z_i}{q^{d_i}-q^{-d_i}}
\lbl{comX_Y},\end{equation}

\be \sum _{r=0}^{1-a_{ij}} (-1)^r \left[ \ba{c} 1-a_{ij} \\ r \ea
\right]_{q^{d_i}} \lambda _{ij}^{a_{ij}+2 r -1 }X_i^{1-a_{ij}-r}
 X_j X_i ^r =0 ,
\lbl{SerreAX}
\end{equation}

\be \sum _{r=0}^{1-a_{ij}} (-1)^r \left[ \ba{c} 1-a_{ij} \\ r \ea
\right]_{q^{d_i}} Y_i^{1-a_{ij}-r} Y_j Y_i ^r =0,
 \lbl{SerreAY} \end{equation} pour $1\leq i,j\leq t $.

Posons
 \be \lbl{comod_AX}  \delta (X_i)= X_i
\otimes 1 +Z_i \otimes E_i ,\end{equation}
\be \lbl{comod_AY}
\delta (Y_i ) = Y_i \otimes K_i^{-1} + 1 \otimes F_i ,\end{equation}
\be \lbl{comod_AZ}
\delta (Z_i ^{\pm 1}) = Z_i ^{\pm 1}\otimes K_i^{\pm 1} ,\end{equation}
pour tout $1\leq i \leq t $.

Nous énonçons maintenant notre résultat principal.

\begin{theo} \lbl{th}
1) Les formules (\ref{comod_AX}), (\ref{comod_AY}) et (\ref{comod_AZ})
munissent $A_{\lambda}$ d'une
structure d'objet galoisien clivé sur $U_q(\mathfrak g)$.

2) Tout objet galoisien sur $U_q(\mathfrak g)$ est homotope à
un objet galoisien de la forme~$A_{\lambda}$.

3) Deux objets $U_q(\mathfrak g)$-galoisiens $A_{\lambda }$ et
$A_{\lambda '}$ sont
homotopes si et seulement si les familles $\lambda $ et
$\lambda '$ les définissant sont égales.
\end{theo}

La suite de l'article est consacrée à la démonstration du
théorème.

\section{Démonstration du théorème}

\subsection{Cocycles sur $\grU$ provenant d'une famille
$\lambda$} Pour toute famille~$\lambda=(\lambda _{ij})_{1\leq i < j
\leq t}$ d'inversibles de $k$, nous notons encore $\lambda
_{ij}=\lambda _{ji}^{-1}$ et~$\lambda _{ii} =1 $ pour tout~$1 \leq
j\leq i \leq t$ et nous définissons un cocycle de~$\grU$ comme suit.
Notons~$\sigma _{\lambda}:k[G]~\times~k[G]~\rightarrow~k$
l'application bilinéaire déterminée par
\be \lbl{defsla1} \sigma_{\lambda} (K_i , K_j) =
\lambda _{ij} , \end{equation}
pour $1\leq i,j \leq t$, et par
\be \lbl{defsla2} \sigma_{\lambda} (g_1g_2,h)= \sigma _{\lambda}
(g_1,h) \sigma _{\lambda} (g_2,h) ,\end{equation} et \be \lbl{defsla3} \sla
(h,g_1g_2) = \sla (h,g_1) \sla (h,g_2) ,\end{equation} pour $g_1,g_2,h \in
G  $. On a donc $\sla (1,g)=\sla (g,1)=1$ pour tout $g\in G$.

Soit $\pi : \grU \rightarrow k[G]$ le
morphisme d'algèbres de Hopf défini par
\be \lbl{def_pi} \pi (E_i)= \pi (F_i)
=0 ,\quad \pi (K_i ^{\pm 1})=K_i^{\pm 1} ,\end{equation} pour $i=1,\ldots ,
t$. Posons~$\widetilde{\sigma _{\lambda}} = \sigma _{\lambda} \circ (\pi
\times \pi): \grU \times  \grU \rightarrow k$. On a le lemme immédiat suivant.
\begin{lemma}
\lbl{lemme_cocycle}
Les applications~$\sigma_{\lambda} $ et~$\widetilde{\sigma _{\lambda}}$ sont des
cocycles normalisés pour les algèbres de Hopf~$k[G]$ et~$\grU$, respectivement, inversibles pour la convolution,
d'inverses respectifs~$\sigma _{\lambda ^{-1}}$ et~$\widetilde{\sigma _{\lambda ^{-1}}}$.
\end{lemma}

\subsection{L'algèbre comodule $A_{\lambda}$ comme $U_q(\mathfrak g)$-extension clivée}
Considérons une famille~$\lambda =(\lambda _{ij}) _{1\leq i,j\leq t}$ d'éléments inversibles de $k$ et les cocycles (pour l'algèbre de Hopf $\grU$)  $\ts$, défini au paragraphe 1.3, et $\rho ^{-1}$ inverse (pour la convolution) du cocycle
$\rho$, défini par Kassel et Schneider et rappelé au paragraphe 1.3. La convolution de ces deux cocycles défini un cocycle~$\sr =\ts * \rho ^{-1}$ pour l'algèbre de Hopf $\grU ^{\rho} \cong \Ug$ (voir \cite[Théorème 5.3]{MS}).

\begin{lemma} Le cocycle inversible $\sr = \widetilde{\sigma _{\lambda}} * \rho ^{-1} : \Ug \otimes \Ug \to k$ vérifie les relations
\be \lbl{normalisation_sr} \sr (1,x)= \varepsilon (x)= \sr (x,1) ,\end{equation}
et si $x,y,z$ appartiennent à la même sous-algèbre $\Ug ^+$ ou
$\Ug ^-$ ou si les termes dans le premier membre de $\sr$ appartiennent à $\Ug ^{-}$ et ceux du second membre appartiennent à $\Ug ^+$, nous avons
\be \lbl{multiplication1_sr} \sr (x,yz) = \sr (x_{(1)},y) \sr (x_{(2)},z), \end{equation}
\be \lbl{multiplication2_sr} \sr (xy,z) = \sr (x,z_{(2)}) \sr (y,z_{(1)}).\end{equation}
De plus, nous avons
\be \lbl{srK_K} \sr (K_i ,K_j )= \lambda _{ij},\end{equation}
\be \lbl{srE_F} \sr (E_i ,F_j ) = \frac{\delta_{ij}}{q^{d_i}-q^{-d_i}},\end{equation}
pour tout $1\leq i,j \leq t$. Pour les autres couples $(x,y)$ de générateurs nous avons
\be \lbl{sr_gene} \sr (x,y)=0.\end{equation}
\end{lemma}

\begin{proof}
La relation de normalisation (\ref{normalisation_sr}) est vérifiée immédiatement à partir des conditions de normalisation pour $\ts$ et $\rho$. Les relations (\ref{multiplication1_sr}) et (\ref{multiplication2_sr}) se déduisent de la définition de $\sigma _{\lambda}$ comme bi-caractère et des relations (\ref{multiplication1_rho}) et (\ref{multiplication2_rho}). Ces relations ne sont donc vérifiées que pour des éléments de la même sous-algèbre~$\Ug ^+$ ou~$\Ug^-$ ou si le premier terme appartient à~$\Ug ^-$ et le second à~$\Ug ^+$ comme pour les relations (\ref{multiplication1_rho}) et (\ref{multiplication2_rho}) de $\rho$.

La relation (\ref{srK_K}) est une conséquence des relations (\ref{comul_U_qK}), (\ref{rhoK_K}) et (\ref{defsla1}).
La relation (\ref{srE_F}) s'obtient à partir des relations (\ref{comul_U_qE}), (\ref{comul_U_qF}), (\ref{rhoE_F}), (\ref{rho_gene}) et (\ref{def_pi}) comme suit :
\be \ba{rcl} \sr (E_i ,F_j ) &= &\sla (E_i ,F_j ) \rho ^{-1} (1, K_j ^{-1}) + \sla (E_i , 1) \rho ^{-1} (1,F_j) \\ &&+ \sla (K_i , F_j) \rho ^{-1} ( E_i ,K_j ^{-1} ) + \sla (K_i ,1)\rho ^{-1}(E_i ,F_j) \\&=& \cfrac{\delta_{ij}}{q^{d_i}-q^{-d_i}},\\ \ea \end{equation}
pour tout $1\leq i,j \leq t$. Pour les autres couples de générateurs le calcul se fait de facon similaire et chaque terme de la somme apparaissant est nul, ce qui implique la relation~(\ref{sr_gene}).
\end{proof}

Posons, pour $1\leq i \leq t$,
\be \lbl{def_varphi} \varphi _{\lambda} (X_i)=E_i,\quad
\varphi _{\lambda} (Y_i)=F_i,\quad
\varphi _{\lambda} (Z_i^{\pm 1} )=K_i^{\pm 1} . \end{equation}

\begin{lemma}
\lbl{lemme_def_phi}
Les formules $(\ref{def_varphi})$ définissent un morphisme
d'algèbres~$U_q(\mathfrak g)$-comodules à droite~$\varphi _{\lambda}: A_{\lambda } \rightarrow\: \Us$.
\end{lemma}

\begin{proof}
a) Vérifions d'abord que $\varphi_{\lambda}$ est un morphisme d'algèbres.
Il suffit d'établir que l'image des relations
$(\ref{comZ_Z}) - (\ref{SerreAY})$ est nulle dans $\Us$.

Considérons la relation (\ref{comZ_Z}) de commutation entre les $Z_i$. On a
\be \lbl{comphiZZ} \ba{rcl}
\varphi _{\lambda} (Z_i Z_j )
&=& \varphi _{\lambda} (Z_i) \ds \varphi _{\lambda} (Z_j) \\
&=& K_i \ds K_j \\
&=& \sr (K_i ,K_j) K_i K_j \\
&=& \lambda _{ij} K_i K_j \\
&=& \lambda _{ij} K_j K_i \\
&=& \lambda _{ij} ^2 \sr (K_j,K_i) K_j K_i \\
&=& \lambda _{ij} ^2 K_j \ds K_i \\
&=& \lambda _{ij} ^2 \varphi _{\lambda} (Z_j) \ds \varphi _{\lambda} (Z_i) \\
&=& \varphi _{\lambda} (\lambda _{ij} ^2 Z_j Z_i ), \\ \ea \end{equation}
pour tout $1\leq i,j \leq t$.
De la m\^eme façon on démontre
\be \varphi _{\lambda} (Z_i Z_i^{-1} )
= \varphi _{\lambda} (Z_i ^{-1} Z_i ) =1 , \end{equation} pour tout $1\leq i,j \leq t$.

Considérons la relation (\ref{comZ_X}) de commutation entre $Z_i$ et $X_j$. On a
\be \ba{rcl}
\varphi _{\lambda} (Z_i X_j )
&=& \varphi _{\lambda} (Z_i) \ds \varphi _{\lambda} (X_j) \\
&=& K_i \ds E_j \\
&=& \sr (K_i ,E_j) K_i 1 + \sr (K_i , K_j) K_i E_j  \\
&=& 0 + \lambda _{ij} K_i E_j \\
&=& \lambda _{ij} q^{d_i a_{ij}} E_j K_i \\
&=& \lambda _{ij} q^{d_i a_{ij}} (0 + \lambda _{ij} \sr (K_j , K_i)  E_j K_i )\\
&=& \lambda _{ij} ^2 q^{d_i a_{ij}} E_j \ds K_i \\
&=& \lambda _{ij} ^2 q^{d_i a_{ij}} \varphi _{\lambda} (X_j) \ds \varphi _{\lambda} (Z_i) \\
&=& \varphi _{\lambda} (\lambda _{ij} ^2 q^{d_i a_{ij}} X_j Z_i )
,\\ \ea \end{equation} pour tout $1\leq i,j \leq t$.

La relation (\ref{comZ_Y}) se démontre de façon similaire.
Pour (\ref{comX_Y}) on a
\be \lbl{comphiXY} \ba{rcl}
\varphi _{\lambda} (X_i Y_j -Y_j X_i  )
&=& \varphi _{\lambda} (X_i) \ds \varphi _{\lambda} (Y_j) - \varphi _{\lambda} (Y_j) \ds \varphi _{\lambda} (X_i)\\
&=& E_i \ds F_j - F_j \ds E_i \\
&=& \left( \sr (E_i ,F_j )1 K_j ^{-1} + \sr (E_i ,1) 1 F_j + \sr (K_i ,F_j) E_i K_j^{-1} \right. \\
&& \left. + \sr (K_i ,1 ) E_i F_j \right)
-\left( \sr (F_j ,E_i ) K_j ^{-1} 1 + \sr (1, E_i )  F_j 1 \right. \\
&& \left. + \sr (F_j,K_i)  K_j^{-1} E_i + \sr (1, K_i  )  F_j E_i \right) \\
&=& \cfrac{\delta _{ij}}{q^{d_i}-q^{-d_i}} K_j^{-1} +  E_i F_j - F_j E_i \\
&=& \delta _{ij} \cfrac{K_i }{q^{d_i} -q^{-d_i} } \\
&=& \varphi _{\lambda} \left( \delta _{ij} \cfrac{Z_i}{q^{d_i} -q^{-d_i} }\right), \\ \ea \end{equation} pour tout $1\leq i,j \leq t$.

Considérons la relation de Serre quantique (\ref{SerreAX}) pour le générateur $X_i$. Notons que
\be \lbl{prodEiEj} \ba{rcl} E_i \ds E_j &=&
\sr (E_i ,E_j) 1 + \sr (E_i ,K_j ) 1 E_j \\
&&+ \sr (K_i , E_j ) E_i 1 + \sr (K_i ,K_j ) E_i E_j \\
&=& \lambda _{ij} E_i E_j ,\\ \ea \end{equation} pour tout $1\leq i,j \leq t$. Nous ne considérons que des éléments de $\Ug ^+$, nous pouvons donc utiliser les relations (\ref{multiplication1_sr}),
(\ref{multiplication2_sr}) et (\ref{sr_gene}) pour calculer les valeurs du cocycle $\sr$, valeurs qui sont nulles dès que le générateur $E_i$ apparait, et obtenir
\be \ba{l} E_i \ds E_j \ds E_j \\\ba{ll} = & \sr ((E_i)_{(1)}, (E_j)_{(1)}) \sr ((E_i)_{(2)} (E_j)_{(2)}, (E_j)_{(1)}) (E_i)_{(3)} (E_j)_{(3)} (E_j)_{(2)} \\=& 0 + \sr (K_i ,K_j ) \sr (K_iK_j, K_j) E_i E_j E_j \\=& \lambda _{ij}^2 E_i E_j ^2 , \ea \ea \end{equation} pour tout $1\leq i,j \leq t$.
De façon similaire, nous avons
\be  E_i \ds E_i \ds E_j = \lambda _{ij}^2 E_i ^2 E_j \end{equation} et
\be E_i \ds E_j \ds E_i = E_i E_j E_i , \end{equation} pour tout $1\leq i,j \leq t$.
Par récurence sur les entiers~$a,b$ et~$c$, on montre facilement que le produit de~$a$ fois le générateur $E_i$, $b$
fois le générateur $E_j$ et $c$ fois le générateur $E_i$ vaut
\be \ba{rcl} \lbl{puissanceE} \underbrace{E_i \ds \cdots \ds E_i}_{a} \ds \underbrace{ E_j \ds \cdots \ds E_j}_{b} \ds \underbrace{E_i \ds \cdots \ds E_i}_{c} &=&E_i^a \ds E_j ^b \ds E_i ^c \\
&=& \lambda _{ij} ^{b(a-c)} E_i ^a E_j ^b E_i ^c  \ea \end{equation} dans~$\Us$ pour tout $1\leq i,j \leq t$. Remarquons en particulier
que la puissance d'un élément $E_i$ est la m\^eme pour le produit
classique de $\Ug$ ou pour celui de $\Us$, ce qui justifie la
notation $E_i ^a$.

Alors, en utilisant (\ref{puissanceE}), on obtient
\be \ba{l} \varphi _{\lambda} \left( \sum _{r=0}^{1-a_{ij}} (-1)^r
\left[ \ba{c} 1-a_{ij} \\ r \ea \right]_{q^{d_i}} \lambda _{ij}^{a_{ij}+2 r-1}  X_i^{1-a_{ij}-r}
 X_j X_i ^r \right)  \\
=  \sum _{r=0}^{1-a_{ij}} (-1)^r \left[ \ba{c} 1-a_{ij} \\ r \ea
\right]_{q^{d_i}}  \lambda _{ij}^{a_{ij}+2 r-1}\varphi
_{\lambda}(X_i)^{1-a_{ij}-r}
\ds  \varphi _{\lambda}(X_j) \ds  \varphi _{\lambda}(X_i) ^r \\
=\sum _{r=0}^{1-a_{ij}} (-1)^r \left[ \ba{c} 1-a_{ij} \\ r \ea
\right]_{q^{d_i}} \lambda _{ij}^{a_{ij}+2 r-1} E_i^{1-a_{ij}-r} \ds
E_j \ds E_i ^r\\
= \sum _{r=0}^{1-a_{ij}} (-1)^r \left[ \ba{c} 1-a_{ij} \\ r \ea
\right]_{q^{d_i}} \lambda _{ij}^{a_{ij} + 2r -1} \lambda _{ij} ^{1-a_{ij}-r -r} E_i ^{1-a_{ij-r} }
E_j E_i ^ r\\
=\sum _{r=0}^{1-a_{ij}} (-1)^r \left[ \ba{c} 1-a_{ij} \\ r \ea
\right]_{q^{d_i}} E_i ^{1-a_{ij}-r} E_j E_i ^ r\\
=0 .
\ea
\end{equation}

De manière similaire, pour les générateurs $F_i$ dans $\Ug ^-$, nous avons
\be \lbl{prodFiFj} \ba{rcl}F_i \ds F_j &=& \sr (F_i ,F_j ) K_i ^{-1} K_j ^{-1} + \sr (F_i , 1)
K_i^{-1} F_j \\&& + \sr (1 , F_j ) F_i K_j^{-1} + \sr (1,1) F_i F_j \\ &=& F_i F_j \\ \ea \end{equation}
et par récurence
\be F_i  ^a\ds F_j ^b \ds F_i ^c = F_i ^a  F_j ^b   F_i ^c \end{equation} pour tout $1\leq i,j \leq t $ et $a,b,c\in \N$.
Le produit des éléments~$F_i$ est donc le m\^eme dans~$\Us$ et dans~$\Ug$. La relation de
Serre quantique (\ref{SerreAY}) dans~$A_{\lambda}$ pour le générateur~$Y_i$ est alors une conséquence immédiate
de la relation~(\ref{SerreF}) pour le générateur~$F_i$ dans~$\Ug$.

b) Pour montrer que $\varphi _{\lambda}$ est aussi un morphisme de comodules, il suffit,
puisque la comultiplication de $\Ug$ est un morphisme d'algèbres, de vérifier
que le diagramme
 \be \ba{ccc} A_{\lambda} &
\xrightarrow{\varphi _{\lambda}} & \Us \\
\downarrow _{\delta} &
& \downarrow _{\Delta } \\
 A_{\lambda } \otimes U_q(\mathfrak g) & \xrightarrow{ \varphi
_{\lambda}  \otimes \rm{Id}}
& \Us \otimes U_q(\mathfrak g) \\
\ea \end{equation} commute pour une famille de générateurs de l'algèbre $A_{\lambda}$.

Pour les générateurs $Z_i ^{\pm 1} $, on a
 \be
\ba{rcl}
  \Delta \circ \varphi _{\lambda } (Z_i ^{\pm 1} )
 &=& \Delta  (K_i ^{\pm 1} ) \\
 &=& K_i ^{\pm 1} \otimes K_i ^{\pm 1} \\
 &=& (\varphi _{\lambda } \otimes \mathrm{Id}) (Z_i ^{\pm 1} \otimes K_i ^{\pm 1}) \\
 &=& (\varphi _{\lambda } \otimes \mathrm{Id}) \circ \delta (Z_i ^{\pm 1}) ;\ea
 \end{equation}
pour les générateurs $X_i$, on a
 \be \ba{rcl} \Delta
\circ \varphi _{\lambda } (X_i )
 &=& \Delta  ( E_i ) \\
 &=& E_i \otimes 1 + K_i \otimes E_i  \\
 &=& (\varphi _{\lambda } \otimes \mathrm{Id}) (X_i \otimes 1 + Z_i \otimes E_i) \\
 &=& (\varphi _{\lambda } \otimes \mathrm{Id}) \circ \delta (X_i );
\ea  \end{equation}
pour les générateurs $Y_i$, on a
 \be \ba{rcl} \Delta
\circ \varphi _{\lambda } (Y_i )
 &=& \Delta  ( F_i ) \\
 &=& F_i \otimes K_i ^{-1}  + 1 \otimes F_i \\
 &=& (\varphi _{\lambda } \otimes \mathrm{Id} )(Y_i \otimes K_i ^{-1} + 1 \otimes F_i) \\
 &=& (\varphi _{\lambda } \otimes \mathrm{Id}) \circ \delta (Y_i )
,\ea \end{equation}
pour tout $1\leq i \leq t$.
\end{proof}

\begin{lemma}\lbl{lemme_phi_iso}
Le morphisme $\varphi _{\lambda} : A_{\lambda } \rightarrow \; \Us$ est un isomorphisme.
\end{lemma}

\begin{proof}
Avec \cite[Chapitre 4]{J}, introduisons l'algèbre de Hopf $U$
engendrée par les générateurs $E_i, F_i$ et $ K_i ^{\pm 1}$ soumis
aux relations (\ref{comK_K}) - (\ref{comE_F}) de commutation de
$U_q(\mathfrak g)$. Cette algèbre peut \^etre munie d'une
structure d'algèbre de Hopf avec la comultiplication, la co\"unité
et l'antipode définies par les m\^emes relations~$(\ref{comul_U_qE}) - (\ref{anti_U_q})$ que pour $\Ug$. L'algèbre
$\Ug$ est alors le quotient de l'algèbre $U$ par l'idéal $I$
engendré par les relations de Serre quantiques~(\ref{SerreE}) et~(\ref{SerreF}) ; notons $P$ le morphisme d'algèbre de Hopf de projection de~$U$ sur~$\Ug$. La famille
\be F_{\underline{i}} ^{\underline{\alpha _i}}
E_{\underline{j}} ^{\underline{\beta_j}} K_{\underline{l} }
^{\underline{\gamma _l}}= F_{i_1}^{\alpha _{i_1}} \cdots
F_{i_n}^{\alpha _{i_n}} E_{j_1}^{\beta _{j_1}} \cdots
 E_{j_p}^{\beta _{j_p}} K_{l_1}^{\gamma _{l_1}}  \cdots
 K_{l_t}^{\gamma _{l_t}} , \end{equation}
o\`u $i_1,\ldots ,i_n , j_1 ,\ldots ,j_p ,l_1,\ldots ,l_t$
parcourent $\{1,\ldots ,t\}$, $\alpha _{i_1} ,\ldots ,\alpha
_{i_n}, \beta _{j_1} ,\ldots ,\beta _{j_p}$ parcourent $\N$ et
$\gamma _{l_1} ,\ldots ,\gamma _{l_t}$ parcourent $\Z$, est une
base de $U$ (voir \cite{J}). L'algèbre de Hopf $k[G]$ est aussi la
sous-algèbre engendrée par les éléments ``group-like'' de~$U$.
Considérons le cocycle (que nous noterons encore $\sr$) $\sr \circ (P\otimes P)
: U\otimes U \to k$ et notons $_{\sr}\! U$ l'algèbre obtenue de~$U$ à partir du cocycle~$\sr$
en suivant le procédé décrit au paragraphe 1.2 et dont le produit est donné par la relation (\ref{produit_tordu_clive_gauche}).

Cherchons maintenant une base de  $_{\sr}\! U$ adaptée au produit~$\ds$.
Remarquons que de la m\^eme manière que pour les relations (\ref{prodEiEj}) et (\ref{prodFiFj}), le produit d'un
élément~$F_i$ avec un élément~$E_j$ vaut dans~$_{\sr}\! U$
\be \lbl{prodEiFj} \ba{rcl}F_i \ds E_j &=& \sr (F_i ,E_j) K_i^{-1} +
\sr (F_i , K_j ) K_i^{-1} E_j \\
&&+ \sr (1 , E_j ) F_i  + \sr (1 ,K_j ) F_i E_j \\&=& F_i E_j ,\\ \ea \end{equation}
celui d'un élément~$E_i$ avec un élément~$K_j$ vaut
\be \lbl{prodFiKj} \ba{rcl} E_i \ds K_j &=& \sr (E_i, K_j)  K_j + \sr (K_i,K_j)  E_i K_j \\ &=&  \lambda _{ij} E_i K_j \ea \end{equation}
et celui de deux éléments~$K_i$ et~$K_j$ vaut\be \lbl{prodKiKj} \ba{rcl}K_i \ds K_j &=& \sr( K_i K_j ) K_i K_j \\&=& \lambda _{ij}
K_i K_j, \\ \ea \end{equation}
pour tout $1\leq i,j\leq t$. Si les indices~$\! j_1 ,\ldots ,j_p ,l_1,\ldots ,l_t$
parcourent $\! \{1,\ldots ,t\}$, notons $J$ la suite d'indice $ j_1,\ldots ,j_n ,l_1 ,$ $l_2, \ldots , l_t$ et
$J_k$ la suite définie à partir de la précédente en ne gardant que les indices à partir
de l'indice~$k$ (par exemple~$J_{j_3}=j_3,\ldots ,j_n ,l_1 ,\ldots , l_t$ et
$J_{l_5}=l_5 ,\ldots , l_t$).

Nous devons calculer les produits dans $_{\sr} \!U$ de produits des générateurs, et nous avons donc besoin des valeurs de $\sr$ sur ces produits. Précisément, nous devons calculer les valeurs du cocycle $\sr$ pour des produits des générateurs~$F_i$, ce qui est possible avec les relations (\ref{multiplication1_sr}) et (\ref{multiplication2_sr}) car nous restons dans la même sous-algèbre $_{\sr} \!U ^-$ (définie de manière évidente comme pour~$\Ug ^-$). De même, la valeur du cocycle $\sr$ sur les produits entre les générateurs $E_j$ et $K_l$ se calcule grâce à ces relations car ces éléments appartiennent à la même sous-algèbre~$_{\sr} \!U ^+$. Enfin, les produits entre les générateurs~$F_i$ et~$E_j$ font intervenir le cocycle $\sigma _{\rho}$ avec comme terme de gauche des éléments de~$_{\sr} \!U ^-$ et comme terme de droite des éléments de $_{\sr} \!U ^+$, ce qui nous permet encore d'utiliser les relations (\ref{multiplication1_sr}) et (\ref{multiplication2_sr}). Alors, de façon similaire à (\ref{prodEiEj}) - (\ref{puissanceE}), nous utilisons les relations~(\ref{prodEiEj}), (\ref{prodFiFj})
et~(\ref{prodEiFj})~-~(\ref{prodKiKj}), exprimant les produits dans~$_{\sr} \!U$ des générateurs~$F_i,E_j$ et~$K_l$ sur la
base~$F_{\ul{i}} ^{\ul{\alpha _i}}  E_{\ul{j}} ^{\ul{\beta _j}}
K_{\ul{l}} ^{\ul{\gamma _l}}$, de manière répétée gr\^ace aux relations~(\ref{multiplication1_sr}) et~(\ref{multiplication2_sr}) de définition de~$\sr$,
pour exprimer les éléments~$F_{\ul{i}} ^{\ul{\alpha _i}} \ds  E_{\ul{j}} ^{\ul{\beta _j}}
\ds K_{\ul{l}}^{\ul{\gamma _l}}$ sur la base~$F_{\ul{i}} ^{\ul{\alpha _i}}  E_{\ul{j}} ^{\ul{\beta _j}} K_{\ul{l}}^{\ul{\gamma _l}}$ :

\be \ba{rcl}
F_{\ul{i}} ^{\ul{\alpha _i}} \ds E_{\ul{j}} ^{\ul{\beta _j}} \ds K_{\ul{l}}^{\ul{\gamma _l}} \!\!
&=& F_{i_1}^{\alpha _{i_1}} \ds \cdots \ds F_{i_n}^{\alpha _{i_n}}\ds E_{j_1}^{\beta _{j_1}} \ds \cdots
\ds E_{j_p}^{\beta _{j_p}}
\\ && \ds K_{l_1}^{\gamma _{l_1}} \ds \cdots \ds K_{l_t}^{\gamma _{l_t}}
\\ &=&  \left(  F_{i_1} ^{\alpha _{i_1}}\cdots   F_{i_n} ^{\alpha _{i_n}} \right) \! \ds \!
\left( (\prod _{j_1\leq j <j' \leq j_p} \beta _{jj'}^{\beta _{j} \beta _{j'}}) E_{j_1}^{\beta _{j_1}} \cdots E_{j_p}^{\beta _{j_p}} \right)
\\ && \ds \left( (\prod _{l_1\leq l <l' \leq l_t} \lambda _{ll'}^{\gamma _l \gamma _{l'}})
K_{l_1}^{\gamma _{l_1}}  \cdots K_{l_t}^{\gamma _{l_t}} \right)
\\   &=&  \prod _{k\in J} \prod _{ k' \in J_k} \lambda _{kk'} ^{\kappa _k \kappa _{k'}}
F_{\ul{i}} ^{\ul{\alpha _i}} E_{\ul{j}} ^{\ul{\beta_j}} K_{\ul{l}}^{\ul{\gamma _l}}, \ea \end{equation}
o\`u $\kappa _k $ désigne~$\beta _k$ si~$k$ est un indice relatif à~$E$ (soit de la forme
$j_m$) et désigne~$\gamma _k$ si~$k$ est un indice relatif à $K$ (soit de la forme $l_m$).
La famille~$E_{\ul{i}} ^{\ul{\alpha _i}} \ds F_{\ul{j}} ^{\ul{\beta _j}} \ds K_{\ul{l}}^{\ul{\gamma _l}}$
forme une base de l'espace vectoriel $U= \:_{\sr}\!U$
puisque les scalaires~$\prod _{k\in J } \prod _{ k'\in J_k} \lambda _{kk'} ^{\kappa _k \kappa _{k'}}$ sont tous non nuls.

Notons~$ Y_{\ul{i}} ^{\ul{\alpha _i}} X_{\ul{j}} ^{\ul{\beta _j}} Z_{\ul{l} }^{\ul{\gamma_l}}$
le produit $Y_{i_1}^{\alpha _{i_1}} \cdots
Y_{i_n}^{\alpha _{i_n}} X_{j_1}^{\beta _{j_1}} \cdots
X_{j_p}^{\beta _{j_p}} Z_{l_1}^{\gamma _{l_1}} \cdots
_{l_t}^{\gamma _{l_t}}$ si $\ul{i},\ul{\alpha _i} ,\ul{j }, $ $ \ul{\beta _j }, \ul{l},\ul{\gamma _l}$
sont les multi-indices correspondant à~$i_1,\ldots ,i_n, \alpha _{i_1} ,\ldots ,\alpha _{i_n},  j_1 ,\ldots , j_p ,$ \\ $ \beta _{j_1},
 \ldots ,\beta _{j_p},l_1,\ldots ,l_t ,\gamma _{l_1},\ldots ,\gamma _{l_t}$ et définissons l'application linéaire~$\psi
: \; _{\sr} U \rightarrow A_{\lambda}$ par sa donnée sur la base
$F_{\ul{i}} ^{\ul{\alpha _i}} \ds  E_{\ul{j}} ^{\ul{\beta _j}}\ds K_{\ul{l}}^{\ul{\gamma _l}}$ :
\be
\psi ( F_{\ul{i}} ^{\ul{\alpha _i}} \ds  E_{\ul{j}} ^{\ul{\beta _j}} \ds K_{\ul{l}}^{\ul{\gamma _l}})
=Y_{\ul{i}} ^{\ul{\alpha _i}} X_{\ul{j}} ^{\ul{\beta _j}} Z_{\ul{l} }^{\ul{\gamma_l}} .\end{equation}

Démontrons que $\psi $ est un morphisme d'algèbres.
Les calculs (\ref{comphiZZ}) - (\ref{comphiXY}) assurent aussi que
\be \lbl{comKKs} K_i \ds K_j = \lambda _{ij}^{2} K_j \ds K_i ,\end{equation}
\be \lbl{comEKs} K_i \ds E_j = \lambda _{ij} ^2 q^{d_i a_{ij}} E_j \ds K_i ,\end{equation}
\be \lbl{comFKs} K_i \ds F_j = q^{-d_i a_{ij}} F_j \ds K_i \end{equation}
et
\be \lbl{comEFs} E_i \ds F_j - F_j \ds E_i = \delta _{ij} \cfrac{K_i }{q^{d_i} -q^{-d_i} }, \end{equation}
pour tout~$1\leq i,j\leq t$.
Alors pour écrire le produit de deux éléments de la base
$F_{\ul{i}} ^{\ul{\alpha _i}} \ds  E_{\ul{j}} ^{\ul{\beta _j}} \ds K_{\ul{l}}^{\ul{\gamma _l}}$ et
$F_{\ul{i'}} ^{\ul{\alpha _i'}} \ds  E_{\ul{j'}} ^{\ul{\beta _j'}} \ds K_{\ul{l'}} ^{\ul{\gamma _l'}}  $, nous utilisons les relations
de commutation~(\ref{comKKs})~-~(\ref{comEFs}) de manière répétée pour obtenir une écriture sur la base :

\be \left( \! F_{\ul{i}} ^{\ul{\alpha _i}} \! \ds \!\!  E_{\ul{j}} ^{\ul{\beta _j}} \! \ds \!\! K_{\ul{l}}^{\ul{\gamma _l}}\! \right) \! \ds \!\!
\left(F_{\ul{i'}} ^{\ul{\alpha _i'}} \! \ds \!\!  E_{\ul{j'}} ^{\ul{\beta _j'}}\! \ds \!\! K_{\ul{l'}} ^{\ul{\gamma _l'}} \right) \!
=\!\! \sum x _{\ul{i''}\ul{j''}\ul{l''}}^{\ul{ \alpha _{i''}}\ul{ \beta _{j''} }\ul{\gamma _{l''}}}
F_{\ul{i''}} ^{\ul{\alpha _{i''}}}\! \ds \!\!  E_{\ul{j''}} ^{\ul{\beta _{j''}}}
\! \ds \!\! K_{\ul{l''}} ^{\ul{\gamma _{l''}}} ,
\end{equation}
avec $x _{\ul{i''}\ul{j''}\ul{l''}}^{\ul{ \alpha _{i''}}\ul{ \beta _{j''} }\ul{\gamma _{l''}}}\in k$.
Remarquons que, dans l'algèbre $A_{\lambda}$, nous avons les mêmes
relations de commutation (\ref{comZ_Z})~-~(\ref{comX_Y}) et donc le produit de deux éléments de la forme
$Y_{\ul{i}} ^{\ul{\alpha _i}}X_{\ul{j}} ^{\ul{\beta _j}} Z_{\ul{l}}^{\ul{\gamma_l}}$ et
$Y_{\ul{i'}} ^{\ul{\alpha_i'}}  X_{\ul{j'}} ^{\ul{\beta _j'}}Z_{\ul{l'}} ^{\ul{\gamma_l'}}$ vaut de la même manière
\be \left(Y_{\ul{i}} ^{\ul{\alpha _i}}X_{\ul{j}} ^{\ul{\beta _j}} Z_{\ul{l}}^{\ul{\gamma_l}}\right)
\left(Y_{\ul{i'}} ^{\ul{\alpha_i'}}  X_{\ul{j'}} ^{\ul{\beta _j'}}Z_{\ul{l'}} ^{\ul{\gamma_l'}} \right) =
\sum x _{\ul{i''}\ul{j''}\ul{l''}}^{\ul{ \alpha _{i''}}\ul{ \beta _{j''} }\ul{\gamma _{l''}}}
Y_{\ul{i''}} ^{\ul{\alpha_i''}}  X_{\ul{j''}} ^{\ul{\beta _j''}}Z_{\ul{l''}} ^{\ul{\gamma_l''}}. \end{equation}   

Nous avons donc
\be \ba{l} \psi ( F_{\ul{i}}^{\ul{\alpha _i}} \ds  E_{\ul{j}} ^{\ul{\beta_j}} \ds K_{\ul{l}} ^{\ul{\gamma _l}} \ds
F_{\ul{i'}} ^{\ul{\alpha _i'}} \ds E_{\ul{j'}} ^{\ul{\beta _j'}} \ds K_{\ul{l'}}^{\ul{\gamma _l'}})\\
\ba{ll}=& \psi \left( \sum x _{\ul{i''}\ul{j''}\ul{l''}}^{\ul{ \alpha _{i''}}\ul{ \beta _{j''} }\ul{\gamma _{l''}}}
F_{\ul{i''}} ^{\ul{\alpha _i''}} \ds  E_{\ul{j''}} ^{\ul{\beta _j''}} \ds K_{\ul{l''}} ^{\ul{\gamma _l''}} \right) \\
=& \sum x _{\ul{i''}\ul{j''}\ul{l''}}^{\ul{ \alpha _{i''}}\ul{ \beta _{j''} }\ul{\gamma _{l''}}}
Y_{\ul{i''}} ^{\ul{\alpha_i''}}  X_{\ul{j''}} ^{\ul{\beta _j''}}Z_{\ul{l''}} ^{\ul{\gamma_l''}} \\
=& \left(Y_{\ul{i}} ^{\ul{\alpha _i}}X_{\ul{j}} ^{\ul{\beta _j}} Z_{\ul{l}}^{\ul{\gamma_l}}
\right)\left(Y_{\ul{i'}} ^{\ul{\alpha_i'}}  X_{\ul{j'}} ^{\ul{\beta _j'}}Z_{\ul{l'}} ^{\ul{\gamma_l'}} \right) \\
=& \psi ( F_{\ul{i}}^{\ul{\alpha _i}} \ds  E_{\ul{j}} ^{\ul{\beta_j}} \ds K_{\ul{l}} ^{\ul{\gamma _l}} )\ds
 \psi (F_{\ul{i'}} ^{\ul{\alpha _i'}} \ds E_{\ul{j'}} ^{\ul{\beta _j'}} \ds K_{\ul{l'}}^{\ul{\gamma _l'}}),\\
\ea \ea \end{equation}
ce qui établit que $\psi $ est un morphisme d'algèbres. Par
suite, comme les générateurs de l'algèbre~$A_{\lambda}$
appartiennent à l'image de $\psi$, celle-ci est surjective. Pour
montrer que $\psi $ se factorise à travers $\Us$, il suffit de
montrer que le noyau de $P$ est inclus dans le noyau de $\psi$.
Soit $u$ appartenant au noyau de $P$ ; alors $u$ appartient à
l'idéal engendré par les relations de Serre quantiques
(\ref{SerreE}) et (\ref{SerreF}). 
Donc $u$ est une combinaison linéaire d'éléments de la forme
\be F_{\ul{i}}^{\ul{\alpha _i}}\! \ds \!\!E_{\ul{j}} ^{\ul{\beta _j}} \! \ds \!\! K_{\ul{l}} ^{\ul{\gamma _l}}
\left(\sum_{r=0}^{1-a_{pq}} (-1)^r \left[ \ba{c} 1-a_{pq} \\ r \ea \right]_{q^{d_p}} E_p^{1-a_{pq}-r}  E_q E_p ^r \right)
F_{\ul{i'}} ^{\ul{\alpha _i'}} \ds E_{\ul{j'}} ^{\ul{\beta _j'}}\! \ds \!\! K_{\ul{l'}}^{\ul{\gamma _l'}} \end{equation}
et de la forme
\be F_{\ul{i}}^{\ul{\alpha _i}}\! \ds \!\!E_{\ul{j}} ^{\ul{\beta _j}}\! \ds \!\! K_{\ul{l}} ^{\ul{\gamma _l}}
\left(\sum \limits _{r=0}^{1-a_{pq}} (-1)^r \left[ \ba{c} 1-a_{pq} \\ r \ea \right]_{q^{d_p}}  F_p^{1-a_{pq}-r} F_p F_q ^r \right)
F_{\ul{i'}} ^{\ul{\alpha _i'}} \! \ds \!\! E_{\ul{j'}} ^{\ul{\beta _j'}}\! \ds \!\! K_{\ul{l'}}^{\ul{\gamma _l'}}, \end{equation}
avec $\ul{i},\ul{\alpha _i},\ul{j},\ul{\beta _j}, \ul{l},
\ul{\gamma _l}$ des multi-indices de la forme précédente et $p,q = 1,\ldots ,t$.
Par conséquent l'élément~$u$ est une combinaison linéaire d'éléments de la forme
\be \ba{l} F_{\ul{i}}^{\ul{\alpha _i}} \! \ds  \!\! E_{\ul{j}} ^{\ul{\beta _j}} \! \ds  \!\! K_{\ul{l}} ^{\ul{\gamma _l}} \! \ds  \!\!
\left( \sum \limits_{r=0}^{1-a_{pq}} (-1)^r \left[ \ba{c} 1-a_{pq} \\ r \ea \right]_{q^{d_p}}\!\!\lambda _{pq}^{a_{pq}+2r -1}
 E_p^{1-a_{pq}-r} \! \ds  \!\! E_q \! \ds  \!\! E_p ^r \right) \\
\ds F_{\ul{i'}} ^{\ul{\alpha _i'}} \ds E_{\ul{j'}} ^{\ul{\beta _j'}}\ds K_{\ul{l'}}^{\ul{\gamma _l'}} \ea \end{equation}
et de la forme
\be \ba{l} F_{\ul{i}}^{\ul{\alpha _i}}\ds E_{\ul{j}} ^{\ul{\beta _j}}\ds K_{\ul{l}} ^{\ul{\gamma _l}} \ds
\left(\sum \limits_{r=0}^{1-a_{pq}} (-1)^r \left[ \ba{c} 1-a_{pq} \\ r \ea \right]_{q^{d_p}}  F_p^{1-a_{pq}-r} \ds F_q \ds F_p ^r \right) \\
\ds F_{\ul{i'}} ^{\ul{\alpha _i'}} \ds E_{\ul{j'}} ^{\ul{\beta _j'}}\ds K_{\ul{l'}}^{\ul{\gamma _l'}}, \ea  \end{equation}
Donc $\psi (u)$ est une combinaison linéaire d'éléments de la forme
\be  \lbl{psiE} \ba{l} \psi \left( F_{\ul{i}} ^{\ul{\alpha _i}} \ds E_{\ul{j}} ^{\ul{\beta _j}} \ds K_{\ul{l}} ^{\ul{\gamma_l}} \ds
\left( \sum \limits_{r=0}^{1-a_{pq}} (-1)^r \left[ \ba{c} 1-a_{pq} \\ r \ea \right]_{q^{d_p}} \lambda _{pq}^{a_{pq}+2r -1} E_p^{1-a_{pq}-r} \ds\right. \right. \\\left. \left.  E_q \ds E_p ^r \right) \ds F_{\ul{i'}} ^{\ul{\alpha '_{i'}}} \ds E_{\ul{j'}} ^{\ul{\beta '_{j'}}} \ds K_{\ul{l'}} ^{\ul{\gamma '_{l'}}} \right) \\
\ba{ll}
=&  \psi \left( F_{\ul{i}} ^{\ul{\alpha _i}} \ds E_{\ul{j}} ^{\ul{\beta _j}} \ds K_{\ul{l}}^{\ul{\gamma_l}}\right)\;
\psi \left(  \sum \limits_{r=0}^{1-a_{pq}} (-1)^r  \left[ \ba{c} 1-a_{pq} \\ r \ea \right]_{q^{d_p}}
 \lambda _{pq}^{a_{pq}+2r -1} \right. \\&\left.  E_p^{1-a_{pq}-r} \ds E_q \ds E_p ^r  \right) \psi \left( F_{\ul{i'}} ^{\ul{\alpha '_{i'}}} \ds E_{\ul{j'}} ^{\ul{\beta '_{j'}}}\ds K_{\ul{l'}} ^{\ul{\gamma '_{l'}}} \right) \ea \\
\ba{ll}
=& Y_{\ul{i}} ^{\ul{\alpha _i}} X_{\ul{j}}^{\ul{\beta _j}} Z_{\ul{l}} ^{\ul{\gamma_l}}
\left(\sum \limits_{r=0}^{1-a_{pq}} (-1)^r\left[ \ba{c} 1-a_{pq} \\ r \ea \right]_{q^{d_p}}
 \lambda _{pq}^{a_{pq}+2r -1}  X_p^{1-a_{pq}-r} X_q X_p ^r \right)
Y_{\ul{i'}} ^{\ul{\alpha '_{i'}}} \\& X_{\ul{j'}} ^{\ul{\beta '_{j'}}} Z_{\ul{l'}} ^{\ul{\gamma '_{l'}}} \ea \ea \end{equation}
et, de façon similaire, d'éléments de la forme
\be \lbl{psiF} \ba{l} \psi(F_{\ul{i}} ^{\ul{\alpha _i}} \ds E_{\ul{j}} ^{\ul{\beta _j}} \ds K_{\ul{l}} ^{\ul{\gamma_l}} \ds
\left(\sum \limits _{r=0}^{1-a_{pq}} (-1)^r \left[ \ba{c} 1-a_{pq} \\ r
\ea \right]_{q^{d_p}}
 F_p^{1-a_{pq}-r} \ds F_q \ds F_p ^r \right)\\
\ds F_{\ul{i'}} ^{\ul{\alpha '_{i'}}} \ds E_{\ul{j'}} ^{\ul{\beta '_{j'}}} \ds K_{\ul{l'}} ^{\ul{\gamma '_{l'}}}) \\
= Y_{\ul{i}} ^{\ul{\alpha _i}} X_{\ul{j}}^{\ul{\beta _j}} Z_{\ul{l}} ^{\ul{\gamma_l}}
\left(\sum \limits _{r=0}^{1-a_{pq}} (-1)^r \left[ \ba{c} 1-a_{pq} \\ r \ea \right]_{q^{d_p}}
 Y_p^{1-a_{pq}-r} Y_q Y_p ^r \right)
Y_{\ul{i'}} ^{\ul{\alpha '_{i'}}} X_{\ul{j'}} ^{\ul{\beta '_{j'}}} Z_{\ul{l'}} ^{\ul{\gamma '_{l'}}} . \ea\end{equation}
Les deux derniers membres de~(\ref{psiE}) et de~(\ref{psiF}) sont nuls en vertu
de~(\ref{SerreAX}) et~(\ref{SerreAY}). Il en résulte que $\psi(u)$ est nul et que le
morphisme $\psi $ se factorise en un morphisme~$\Psi : \: \Us \rightarrow A_{\lambda}$.

Le morphisme $\Psi$ est surjectif puisque $\psi$ l'est. De plus nous avons la relation
\be \ba{rcl}
(\varphi _{\lambda} \circ \Psi) (
F_{\ul{i}} ^{\ul{\alpha _i}} \ds E_{\ul{j}} ^{\ul{\beta _j}} \ds K_{\ul{l}} ^{\ul{\gamma_l}} )& = &
\varphi _{\lambda} (Y_{\ul{i}} ^{\ul{\alpha _i}}
X_{\ul{j}} ^{\ul{\beta _j}} Z_{\ul{l}} ^{\ul{\gamma_l}} ) \\
&=& \varphi _{\lambda
}(Y_{\ul{i}}) ^{\ul{\alpha _i}}
\ds \varphi _{\lambda }( X_{\ul{j }})^{\ul{\beta _j}} \ds
\varphi _{\lambda }( Z_{\ul{l}} ) ^{\ul{\gamma_l}} \\
&=& F_{\ul{i}} ^{\ul{\alpha _i}} \ds E_{\ul{j}} ^{\ul{\beta _j}}
 \ds K_{\ul{l}} ^{\ul{\gamma_l}}, \ea \end{equation}
pour tout multi-indices $\ul{i},\ul{\alpha _i},\ul{j},\ul{\beta _j},\ul{j} ,\ul{\gamma _l}$.
Donc l'application~$\varphi _{\lambda} \circ \Psi$ est égale à l'identité de~$\Us$.
En conséquence, $\Psi$ est injective. Comme $\Psi$ est surjective, elle est bijective d'inverse~$\varphi _{\lambda}$.
\end{proof}

\subsection{Démonstration du théorème}

Le point $(1)$ est une conséquence du lemme~\ref{lemme_phi_iso}.

Démontrons le point $(2)$. Soit $A$ un objet galoisien de $\Ug$.
Notons $i$ le plongement naturel de~$k[G]$ dans~$\Ug$ ; ce morphisme d'algèbres de Hopf nous
permet en particulier de définir une structure d'algèbre~$\Ug$-comodule
à gauche sur $k[G]$ par 

\begin{equation}  k[G] \xrightarrow{\Delta _{k[G]}} k[G]\otimes k[G] \xrightarrow{i \otimes \rm{Id}} U_q(\mathfrak g) \otimes k[G] . \end{equation} 

Comme nous l'avons expliqué au paragraphe
1.1, ce plongement $i : k[G] \to \Ug$ induit une application 

\begin{equation}  i^{\star} : \mathrm{Gal} _k(U_q(\mathfrak g)) \rightarrow \mathrm{Gal} _k (k[G]). \end{equation}

Considérons l'image $i^{\star}(A)$ de $A$ dans $\mathrm{Gal} _k
(k[G])$. D'après \cite[Prop. 3.2]{KS},
\begin{equation} \lbl{galhom} \Gal (k[G]) \cong \mathrm{H} ^2 (G,k^{\star}).\end{equation}
Il est bien connu que ce dernier groupe est isomorphe \`a $\mathrm{Hom}(\Lambda^2 \Z ^t ,k^*)$ (voir par exemple
\cite[Théorème~V.6.4~(iii)]{B}). Par cons\'equent, il existe une famille~$\lambda$ telle que~$i^*(A)
\cong \: _{\sigma_{\lambda}}\! k[G]$. 

Le lemme \ref{lemme_iUs} assure que $i^*(A) \cong \: _{\sigma_{\lambda}}\! k[G]$ est isomorphe, comme algèbre, à~$_{\sigma _{\rho}} U_q(\mathfrak g) \Box _{U_q(\mathfrak g)} k[G]$ qui vaut $i^*(_{\sigma _{\rho}} U_q(\mathfrak g))$, par définition de $i^*$. 
Les lemmes~\ref{lemme_def_phi} et~\ref{lemme_phi_iso} assurent que $_{\sigma _{\rho}} U_q(\mathfrak g)  \cong A_{\lambda}$ et par suite nous avons  
\begin{equation} 
i^{\star}(A) \cong i^{\star} (A_{\lambda}). \end{equation}
Or Kassel et Schneider \cite[Théorème 4.5]{KS} ont montré que
l'application~$i^{\star}: \cH _k(\Ug ) \rightarrow \cH _k (k[G])$
est une bijection sur les classes d'homotopie d'extensions
galoisiennes. Par conséquent, $A$ et $A_{\lambda}$ sont homotopes.

Démontrons le point $(3)$. Supposons que $A_{\lambda} $ et
$A_{\lambda '}$ définissent le m\^eme élément de $\cH (\Ug )$.
Alors $i^{\star} (A_{\lambda }) = i^{\star } (A_{\lambda '}) $
dans $\cH _k(k[G]) $. D'après \cite[Prop 3.2]{KS}, on a $\cH
_k(k[G]) \cong \Gal _k(k[G])$. Il résulte de ceci et de la
bijection (\ref{galhom}) que~$\lambda =\lambda '$.

\backmatter


\begin{thebibliography}{EM}
\rm
\bibitem[B]{B} K. S. Brown, Cohomology of groups, Graduate Texts in Mathematics, volume~87,
Springer-Verlag, New York-Berlin, 1982.

\bibitem[EM]{EM} S. Eilenberg, J.-C. Moore,
Homology and fibrations. I. Coalgebras, cotensor product and its derived functors,
Comment. Math. Helv. {\bf 40} (1966) 199--236.
\bibitem[J]{J}
J.-C. Jantzen, Lectures on Quantum Groups,
Graduate Studies in Mathematics, Volume 6, Amer. Math. Soc.,
Providence, RI, (1995).

\bibitem[K]{K}
C. Kassel, Quantum principal bundles up to homotopy equivalence,
The Legacy of Niels Henrik Abel, The Abel Bicentennial, Oslo, 2002,
O.~A.~Laudal, R.~Piene (eds.), Springer-Verlag (2004), 737--748.


\bibitem[KS]{KS} C. Kassel, H.-J. Schneider, Homotopy theory of Hopf Galois
extensions. {\it Annales Inst. Fourier} (Grenoble), {\bf 55} (2005), 2521--2550.


\bibitem[M]{M}
S. Montgomery, Hopf algebras and their actions on rings,
CBMS Regional Conference Series in Mathematics, 82,
American Mathematical Society, Providence, RI, (1993).

\bibitem[MS]{MS}
S. Montgomery, H.-J. Schneider, Krull relations in Hopf Galois extensions : lifting and twisting, J.~Algebra {\bf 288} (2005), 364--383.

\bibitem[S]{S}
H.-J. Schneider, Principal homogeneous spaces for arbitrary Hopf algebras,
Israel J. Math.  {\bf 72}  (1990),  no. 1-2, 167--195.

\end{thebibliography}
\end{document}